\documentclass{article}

\usepackage{amsmath}
\usepackage{url}

\begin{document}

\title{How would Riemann evaluate $\zeta(2n)$?}
\author{Marco Dalai}
\date{}
\maketitle

Let $\zeta(s)$ be the Riemann zeta function defined, for $\Re(s)>1$, by the sum
\begin{equation}\label{zeta}
  \zeta(s):=\sum_{k=1}^\infty\frac{1}{k^s}.
\end{equation}

In his famous paper \cite{riemann}, Riemann introduces the integral representation 
\begin{equation}\label{ident:intmenozeta}
\zeta(s) = \frac{1}{\Gamma(s)}\int_0^\infty
\frac{x^{s-1}}{e^x-1}\,dx, \qquad \Re(s)>1,
\end{equation}
and then, building upon this representation, he defines $\zeta(s)$ for all complex values of $s$ by using the contour integral equation
\begin{equation}\label{ident:intcontour}
2 \sin(2 \pi s)\Gamma(s) \zeta(s) = i \int_\infty^\infty
\frac{(-x)^{s-1}}{e^x-1}\,dx.
\end{equation}

A remarkable well-known result, first obtained by Euler in 1735, is the evaluation of $\zeta(s)$ for even integer values of $s$, that is
\begin{equation}\label{euler}
\zeta(2n)=\frac{(2\pi)^{2n}(-1)^{n+1} B_{2n}}{2(2n)!},
\end{equation}
where $B_{i}$ is the $i$-th Bernoulli number. A special case of particular historical importance is the relation $\zeta(2)=\pi^2/6$, for which today several proofs are known. In order not to overload the bibliography, we refer to \cite{deamo} for a recent new proof and for useful references of previous ones. We also mention that fourteen different proofs for this special case have been collected by R. Chapman in an unpublished note \cite{chapman}.

In his book on the Riemann's zeta function \cite[p. 12]{edwards}, Prof. Edwards observes that {\em ``There is no easy way to deduce this famous formula of Euler's}  [eq. \eqref{euler} here] {\em from Riemann's integral formula for $\zeta(s)$} [equation \eqref{ident:intcontour} here] {\em and it may well have been this problem of deriving} [eq. \eqref{euler} here] {\em anew which led Riemann to the discovery of the functional equation of the zeta function}''.
It is not unfair to conjecture about possible attempts by Riemann to find new ways to evaluate $\zeta(2n)$, since most of the currently known methods became popular much later.
The aim of this note is to show that there is indeed an easy way to deduce \eqref{euler} directly from \eqref{ident:intmenozeta}, a fact that appears to be interesting, at least from a historical point of view, but doesn't seem to have been previously observed to the best of the author's knowledge.

Note that, even if this proof relies on contour integration in the complex domain, no use is made of the residue theorem and only Cauchy's theorem is used. Thus, even in the special case $s=2$, this proof is different from all those presented in \cite{chapman}.

First note that, by a change of variable $x=2\tau$ in the integral appearing on the right hand side of \eqref{ident:intmenozeta}, and then using the identity 
\begin{equation}
\frac{2}{e^{2\tau}-1}=\frac{1}{e^{\tau}-1}-\frac{1}{e^\tau+1},
\end{equation}
we obtain, after simple manipulations,
\begin{equation}
\int_0^\infty \frac{\tau^{s-1}}{e^\tau+1}\,d\tau = (1-2^{1-s})\Gamma(s)\zeta(s).
\label{ident:intpiuzeta}
\end{equation}

Now, for complex $z=x+iy$ and integer $s$, let $f(z)=z^{s-1}/(e^z-1)$ and, given $R>0$, let
$\Gamma_R$ be the rectangular contour of vertices $0$, $R$,
$R+i\pi$ and $i\pi$ in the complex plane. By Cauchy's theorem, we
have that
\begin{equation}
\int_{\Gamma_R}f(z)dz=0\label{inteq0}.
\end{equation}
It is easy to see that the integral over the right vertical side
of the contour tends to zero when $R\to \infty$. So, when $R\to
\infty$, equation (\ref{inteq0}) gives
\begin{equation}
 \int_0^\infty
 \frac{x^{s-1}}{e^x-1}\,dx-\int_0^\infty\frac{(x+i\pi)^{s-1}}{e^{x+i\pi}-1}\,dx=
 i\int_0^\pi\frac{(iy)^{s-1}}{e^{iy}-1}\,dy.
 \label{eq:genric}
\end{equation}
For the particular case $s=2$, using \eqref{ident:intmenozeta} and \eqref{ident:intpiuzeta}, this equation reads
\begin{equation}
  \zeta(2)+\frac{1}{2}\,\zeta(2)+i\pi\int_0^\infty\frac{1}{e^x+1}\,dx=
  \int_0^\pi\frac{y}{2}\,dy +\frac{i}{2}\int_0^\pi \frac{y
  \sin y}{1-\cos y}\,dy.\nonumber
\end{equation}
Taking the real part of both sides one finds the desired result $\zeta(2)=\pi^2/6$.
Similarly, in the case of even integer $s=2n$, using \eqref{ident:intmenozeta} and \eqref{ident:intpiuzeta}, after some manipulations the real part of \eqref{eq:genric} gives

\begin{equation}
\Gamma(2n)\zeta(2n)+\sum_{k=0}^{n-1} \alpha(n,k) \zeta(2n-2k)=(-1)^{n-1}\frac{\pi^{2n}}{4n}
\end{equation}
where
\begin{equation}
\alpha(n,k)=\left(1-2^{1-2n+2k}\right)(-\pi^2)^k\binom{2n-1}{2k}\Gamma(2n-2k).
\end{equation}
This allows us to obtain the values $\zeta(2n)$ recursively, the solution being then expressible in closed form with the use of Bernoulli numbers.

It is reasonable to speculate that Riemann could have noticed this procedure while working out the analytic extension of $\zeta(s)$ by means of \eqref{ident:intcontour} using \eqref{ident:intmenozeta} as a starting point.
It is worth pointing out that, as a side result, taking the imaginary part of \eqref{eq:genric} for even $s$ or the real part for odd $s$, one obtains integral and series representations for the values $\zeta(2k+1)$. This also suggests an intuitive interpretation of the different nature of $\zeta(2n)$ and $\zeta(2n+1)$, being the quantities associated respectively with the imaginary and the real parts of the function $f(z)=z/(e^z-1)$ on the imaginary axis.


\bigskip

\noindent\textit{Department of Information Engineering, 
University of Brescia, Italy\\
marco.dalai@ing.unibs.it}


\begin{thebibliography}{9}

\bibitem{chapman} R. Chapman, Evaluating $\zeta(2)$ (1999), available at\\ 
\url{http://empslocal.ex.ac.uk/people/staff/rjchapma/etc/zeta2.pdf}

\bibitem{deamo}
E. De Amo, M. D\'iaz Carrillo and J. Fern\'andez-S\'anchez, Another proof of Euler's formula for $\zeta(2k)$, \emph{Proc. Amer. Math. Soc.}, \textbf{139} (2011) 1441-1444.

\bibitem{edwards}
H. M. Edwards,\textit {Riemann's Zeta Function}. Academic Press, New York, 1974.


\bibitem{riemann}
G. F. B. Riemann, \"{U}ber die anzahl der primzahlen unter einer gegebenen gr\"{o}sse, \emph{Monatsber. Berlin. Akad.}, (1859) 671-680.

%
%
\end{thebibliography}
\end{document}